\documentclass{ifacconf}

\usepackage{graphicx}      
\usepackage{natbib}        

\input{custom.sty}
\input{ifacthm.sty}

\newcommand{\xd}{x_{[d]}}
\newcommand{\xdt}{x_{[d]}^\Tr}
\newcommand{\pcite}[1]{(\cite{#1})}
\newcommand{\setS}{\mathcal{S}}
\newcommand{\piB}{\pi_{\Image(B)^\perp}}

\newcommand{\polyset}[1]{
  \begin{subfigure}[b]{0.22\textwidth}
    \centering
      \includegraphics[trim=3.7cm 0cm 3.5cm 0cm, clip, width=\textwidth]{primal_#1.png}
  \end{subfigure}%
  ~%
  \begin{subfigure}[b]{0.22\textwidth}
    \centering
      \includegraphics[trim=3.7cm 0cm 3.5cm 0cm, clip, width=\textwidth]{polar_#1.png}
  \end{subfigure}
}

\begin{document}

\begin{frontmatter}

\title{Geometric control of algebraic systems
}

\thanks[footnoteinfo]{RJ is a FNRS honorary Research Associate. This project has received funding from the European Research Council (ERC) under the European Union's Horizon 2020 research and innovation programme under grant agreement No 864017 - L2C. RJ is also supported by the Walloon Region, the Innoviris Foundation, and the FNRS (Chist-Era Druid-net).}

\author[First]{Beno\^it Legat}
\author[Second]{Rapha\"el M. Jungers}

\newcommand{\email}[1]{\emph{\texttt{#1}}}
\address[First]{ICTEAM institute, UCLouvain, Louvain-la-Neuve, Belgium. (e-mail: \email{benoit.legat@uclouvain.be}).}
\address[Second]{ICTEAM institute, UCLouvain, Louvain-la-Neuve, Belgium. (e-mail: \email{raphael.jungers@uclouvain.be}).}

\begin{abstract}
  \,\,\,\,In this paper, we present a geometric approach for computing the
  controlled invariant set of a continuous-time control system.
  While the problem is well studied for in the ellipsoidal case,
  this family is quite conservative for constrained or switched linear systems.
  We reformulate the invariance of a set as an inequality for its support function
  that is valid for any convex set.
  This produces novel algebraic conditions for the invariance of
  sets with polynomial or piecewise quadratic support function.
  We compare it with the common algebraic approach for polynomial sublevel sets
  and show that it is significantly more conservative than our method.
\end{abstract}


\end{frontmatter}

\section{Introduction}

Computing\footnote{A preliminary version of this work appears in \cite{legat2020set}.} controlled invariant set is paramount in many applications~\pcite{blanchini2015set}.
Indeed, the existence of a controlled invariant set is equivalent to the stabilizability of a control system~\pcite{sontag1983lyapunov} and
a (possibly nonlinear) stabilizable state feedback can be deduced from the controlled invariant set~\pcite{barmish1985necessary}.

The stabilizability of a linear time-invariant (LTI) control system is equivalent
to the stability of its uncontrollable subspace (which is readily accessible
in its Controllability Form) \cite[Section~2.4]{wonham1974linear}.
Indeed, the eigenvalues of its controllable subspace can be fixed to any value
by a proper choice of linear state feedback.
The resulting controlled system is stable hence an invariant ellipsoid can be
determined by solving a system of linear equations~\pcite{Liapounoff1907}.
This set is also controlled invariant for the control system.
When a control system admits an ellipsoidal controlled invariant set,
it is said to be \emph{quadratically stabilizable}.
When there exists a linear state feedback such that the resulting autonomous
system admits an ellipsoidal invariant set, it is said to be
\emph{quadratically stabilizable via linear control}.

While the stabilizability of LTI control systems is equivalent to their
quadratic stabilizability via linear control,
it is no longer the case for \emph{uncertain} or \emph{switched} systems \pcite{petersen1985quadratic}.
Furthermore, it is often desirable for \emph{constrained} systems to find a controlled invariant set of maximal
volume (or which is maximal in some direction~\pcite{ahmadi2018robust}).
For such problem, the method detailed above is not suitable as it does not
take any volume consideration but more importantly, the maximal volume
invariant set may not be an ellipsoid and may not be rendered stable via
a linear control.
For this reason, the Linear Matrix Inequality (LMI) \eqref{eq:controlled_LMI} was devised to encapsulate the
controlled invariance of an ellipsoid via linear control \cite[Section~7.2.2]{boyd1994linear}
and the conservatism of the choice of linear control was analysed \pcite{sontag1983lyapunov}.
As the linearity of the control was found to be conservative for uncertain systems \pcite{petersen1985quadratic},
the LMI \eqref{eq:invLMI} was found to encapsulate controlled invariance of an ellipsoid
via \emph{any} state-feedback \cite{barmish1985necessary}.

These LMIs have had a tremendous impact on control, but unfortunately the approach is limited to ellipsoids due to its algebraic nature.
We reinterpret it in a geometric/behavioural framework,
based on convex analysis,
which allows us to formulate a general condition for the controlled invariance of arbitrary
convex sets via any state-feedback in \theoref{invconvex}.
While this condition reduces to \eqref{eq:invLMI} for the special case of ellipsoids,
it provides a new methods for computing controlled invariant convex sets with convex polynomial and piecewise quadratic support functions.

In \secref{alg}, we review the classical LMI's for the invariance and controlled invariance of ellipsoids and discusses the challenges for its generalization to sublevel sets of polynomials.
In \secref{geometric}, we develop a generic condition of control invariance
for continuous-time systems using our geometric approach.
We particularize it for ellipsoids (resp. sets with polynomial and piecewise quadratic support functions)
in \secref{ell} (resp. \secref{poly} and \secref{pell}).
We illustrate these new results with a numerical example in \secref{exem}.

\paragraph*{Reproducibility}
The code used to obtain the results is
published on codeocean \pcite{legat2021continuousCO}.
The set programs are reformulated
by
SetProg \pcite{legat2019set}
into
a Sum-of-Squares program
which is
reformulated into a semidefinite program
by SumOfSquares \pcite{weisser2019polynomial}
which is
solved by Mosek v8 \pcite{mosek2017mosek81034} through MathOptInterface \pcite{legat2020mathoptinterface}.


\section{Algebraic approach}
\label{sec:alg}

Computing an ellipsoidal invariant set for an autonomous system
$\dot{x} = Ax$ where $A \in \R^{n \times n}$ consists in searching for an ellipsoidal set
\[ \mathcal{E}_P = \{\, x \in \R^n \mid x^\Tr P x \leq 1 \,\} \]
that satisfies the \emph{Nagumo condition} \cite[Theorem~4.7]{blanchini2015set}: 
$x^\Tr PAx \leq 0$ for any $x \in \mathbb{R}^n$.
The Nagumo condition for ellipsoids is equivalent to the
LMI:
\begin{equation}
  \label{eq:uncontrolled_LMI}
  A^\Tr P + PA \preceq 0.
\end{equation}
which allows to search for ellipsoidal invariant sets using
semidefinite programming.

Consider the continuous-time control linear system
\begin{equation}
  \label{eq:cs}
  \dot{x} = Ax + Bu
\end{equation}
where $A \in \R^{n_x \times n_x}$ and $B \in \R^{n_x \times n_u}$
with the following definition of invariance.

\begin{mydef}[Controlled invariant set]
  \label{def:cis}
  A set $S$ is \emph{controlled invariant} for system \eqref{eq:cs}
  if for any state $x_0 \in S$, there exists a control $u(t)$ such that
  the trajectory with initial state and control $u$ remains in $S$.
\end{mydef}

With the presence of the control $u$ in the system,
the Nagumo condition becomes:
\begin{equation}
  \label{eq:existential}
  \forall x \in \R^{n_x}, \exists u \in \R^{n_u}, x^\Tr P(Ax + Bu) \leq 0.
\end{equation}
The control term $u$, or more precisely the existential quantifier $\exists$ prevents
us to transform this into an LMI directly.

Fixing the control to a linear state feedback $u(x) = Kx$ for some matrix $K$ allows
to fallback to the case of autonomous system $\dot{x} = (A + BK) x$.
Then, the invariance condition can be formulated
as the \emph{Bilinear Matrix Inequality} (BMI):
\[ A^\Tr P + PA + K^\Tr B^\Tr P + PBK \preceq 0 \]

While the matrix inequality is bilinear in $K$ and $P$, and BMI's are NP-hard to solve in general \pcite{toker1995np},
a clever algebraic manipulation
allows to reformulate it as a
Linear Matrix Inequality (LMI) in $Q := P^{-1}$ and $Y := KQ$,
where the sought controlled invariant ellipsoid is given by $\mathcal{E}_P$,
see e.g. \cite[Section~7.2.1, Section~7.2.2]{boyd1994linear} and \cite[Section~4.4.1]{blanchini2015set} for more details.
The linear matrix inequality is
\begin{equation}
  \label{eq:controlled_LMI}
  QA^\Tr + AQ + Y^\Tr B^\Tr + BY \preceq 0.
\end{equation}

Because the algebraic manipulation which allows to reformulate the BMI into
an LMI is done at the level of matrices, it is not clear how this
approach can be generalized to other families of sets such as the
ones considered in \secref{poly} and \secref{pell}.
Moreover, searching for ellipsoidal controlled invariant sets may be rather
restrictive and the conservativeness is amplified for the class of hybrid systems.

One attempt to generalize it to controlled invariant sublevel sets of polynomials of degree $d$
is developed in \cite{Prajna2004}.
While the method allows to consider systems defined by polynomial equations,
we show below that for linear systems, it
has significant restrictions for $d > 2$
and the case $d = 2$ reduces to the ellipsoidal case given by \eqref{eq:controlled_LMI}.
This suggests that the methods is essentially restricted to compute invariant sublevel sets of polynomials
of a most twice the degree of the polynomial equations.
On the other hand, the condition we find in \cororef{poly} is necessary and sufficient.
The remaining of this section particularize the approach for linear systems,
more details can be found in \cite{Prajna2004}.


Let $\xd$ denote the vector of all monomials of degree $d$ with the variables $x_i$.
Consider the set $\setbuild{x}{\xdt P \xd \le 1}$ for a symmetric positive definite matrix $P$
and a state feedback of the form
$u(x) = K_1(x)x + K_2(x)\xd$
where $K_1(x), K_2(x)$ are matrices of the appropriate dimensions whose
entries are polynomials in $x$.
The invariance of the set for the autonomous system
$\dot{x} = Ax + BK_1(x)x + BK_2(x)\xd$
is equivalent to the nonnegativity of the polynomial
\[
  \xdt P M(x) B K_2(x) \xd + \xdt PF(x) x
\]
where $F(x) = M(x)BK_1(x) + M(x)A$
and $M(x)$ is the jacobian of the transformation $x \mapsto \xd$.
This polynomial can be rewritten in the matrix form:
\begin{equation}
  \label{eq:polyneg}
  \begin{bmatrix}
    \xd\\
    x
  \end{bmatrix}^\Tr
  \begin{bmatrix}
    P M(x) B K_2(x) & PF(x)\\
    0 & 0
  \end{bmatrix}
  \begin{bmatrix}
    \xd\\
    x
  \end{bmatrix}.
\end{equation}
The following is therefore a sufficient condition for the invariance:
\begin{equation}
  \label{eq:polybmi}
  \forall x \in \R^n, \quad
  \begin{bmatrix}
    2P M(x) B K_2(x) & PF(x)\\
    F^\Tr(x)P & 0
  \end{bmatrix} \preceq 0.
\end{equation}
After a similarity transformation with the block diagonal matrix
$\BlockDiag(P^{-1}, I_{n_x})$ where $I_{n_x} \in \R^{n_x \times n_x}$ is the identity matrix, the condition is rewritten as:
\begin{equation}
  \label{eq:polyLMI}
  \forall x \in \R^n, \quad
  \begin{bmatrix}
    2M(x) B Y(x) & F(x)\\
    F^\Tr(x) & 0
  \end{bmatrix} \preceq 0
\end{equation}
where $Y(x) = K_2(x)P^{-1}$.

By \propref{block}, for \eqref{eq:polyLMI} to hold we need all entries of $F(x)$
to be zero polynomials which is quite restictive.
The example below discusses the
conservativeness of \eqref{eq:polyneg} compared to \eqref{eq:polybmi}.

\begin{myexem}
  Consider the autonomous system $\dot{x} = -x$ and the invariant set $\setbuild{x}{x^4 \le 1}$.
  The condition \eqref{eq:polyneg} is satisfied:
  \[
    \begin{bmatrix}
      x^2\\
      x
    \end{bmatrix}^\Tr
    \begin{bmatrix}
      0 & -2x\\
      0 & 0
    \end{bmatrix}
    \begin{bmatrix}
      x^2\\
      x
    \end{bmatrix}
    =
    -x^4
  \]
  while the matrix
  \[
    \begin{bmatrix}
      0 & -2x\\
      -2x & 0
    \end{bmatrix}
  \]
  of conditions \eqref{eq:polybmi} and \eqref{eq:polyLMI} is indefinite for any nonzero $x$.
\end{myexem}

\section{Geometric approach}
\label{sec:geometric}

In this section we derive a characterization of
the controlled invariance of a closed convex set
under the form of an inequality for its support function.
We start by showing the equivalence of the notion of invariance with
another class of systems that directly models
the geometric behaviours of the trajectories of
control systems with unconstrained input.

Consider the continuous-time \emph{algebraic} linear systems:
\begin{equation}
  \label{eq:as}
  E\dot{x} = Cx.
\end{equation}
with the following definition of invariance.

\begin{mydef}[Invariant set for an algebraic system]
  \label{def:is}
  A set $S$ is \emph{invariant} for system \eqref{eq:as}
  if for any state $x_0 \in S$, there exists a trajectory of the system that
  remains in $S$.
\end{mydef}

Note the use of ``there exists'' instead of ``for all'' in the definition of invariance
as both versions exists; see \cite[Remark~4]{legat2020sum} for more details.

\begin{myprop}
  \label{prop:projB}
  Given a subset $\mathcal{S} \subseteq \R^n$ and
  matrices $A \in \R^{r \times n}, B \in \R^{r \times m}$, the following holds:
  \[ A \mathcal{S} + B \R^m = \piB{}^{-1} \piB{} A \mathcal{S} \]
  where $\piB{}$ is any projection matrix onto the orthogonal subspace of $\Image(B)$.
  \begin{proof}
    Given $x \in \mathcal{S}$ and $y \in \R^r$, we have $y \in A \{x\} + B \R^m$ if and only if $y - Ax \in \Image(B)$ or equivalently, $\piB y = \piB Ax$.
  \end{proof}
\end{myprop}

\begin{myprop}
  \label{prop:csas}
  Consider a control system~\eqref{eq:cs}
  and an arbitrary projection matrix $\piB{}$ onto the orthogonal subspace of $\Image(B)$.
  A set $S$ is controlled invariant for the control system~\eqref{eq:cs} with $\mathcal{U} = \R^{n_u}$,
  as defined in \defref{cis},
  if and only if
  it is controlled invariant for the algebraic system 
  \[ \piB{}\dot{x} = \piB{}Ax, \]
  as defined in \defref{is}.
  \begin{proof}
    By \propref{projB}, there exists $u \in \R^{n_u}$ such that $\dot{x} = Ax + Bu$
    if and only if $\piB{}\dot{x} = \piB{}Ax$.
    As the input $u$ is unconstrained, the result follows.
  \end{proof}
\end{myprop}

The Nagumo condition 
for algebraic systems has the following form.

\begin{myprop}
  \label{prop:nagumo}
  A closed set $\setS$ is invariant for system~\eqref{eq:as},
  as defined in \defref{is},
  if and only if
  \begin{equation}
    \label{eq:invtang}
    \forall x \in \partial \setS, \exists y \in T_S(x), Ey = Ax.
  \end{equation}
\end{myprop}





See \secref{convex} for a brief review of the concepts
of convex geometry needed for the remaining of this section.
The invariance condition \eqref{eq:invtang} can be rewritten in terms
of \emph{exposed faces}.

\begin{mytheo}[Controlled invariance of convex set]
  \label{theo:invconvex}
  A convex set $\Csetvar$ is invariant for system \eqref{eq:as}
  with matrices $C, E \in \mathbb{R}^{r \times n}$
  if and only if
  \begin{equation}
    \label{eq:invface}
    \forall z \in \mathbb{R}^r, \forall x \in F_\Csetvar(E^\Tr z), \la z, Cx \ra \le 0.
  \end{equation}
  \begin{proof}
    As $\Csetvar$ is convex, $T_\Csetvar(x)$ is a convex cone. By definition
    of the polar of a cone, $x \in ET_\Csetvar(x)$ if and only if
    $\langle y, x \rangle \le 0$ for all $y \in \polar{[ET_\Csetvar(x)]}$.
    By \propref{podu}, $\polar{[ET_\Csetvar(x)]} = E^{-\Tr}N_\Csetvar(x)$.
    Therefore, the set $\Csetvar$ is invariant if and only if
    \begin{equation}
      \label{eq:invnorm}
      \forall x \in \partial \Csetvar, \forall z \in E^{-\Tr} N_\Csetvar(x), \langle z, Cx \rangle \le 0.
    \end{equation}
    By \propref{normalsupfun}, we have
    \begin{multline*}
      \{\, (x, z) \in \partial \Csetvar \times \mathbb{R}^r \mid E^\Tr z \in N_\Csetvar(x) \,\} = \\
      \{\, (x, z) \in \partial \Csetvar \times \mathbb{R}^r \mid x \in F_\Csetvar(E^\Tr z) \,\}.
    \end{multline*}
  \end{proof}
\end{mytheo}

As we show in the remaining of this section, \theoref{invconvex} allows to reformulate
the invariance as an inequality in terms of the support function.
This allows to combine the invariance constraint with other set constraints that can be formulated in terms of support functions.
Moreover, for an appropriate family of sets, also called \emph{template}, the set program
can be automatically rewritten into a convex program combining all constraints using \emph{set programming} \cite{legat2019set, legat2020set}.
For this reason, we only focus on the invariance constraint and do not detail how to
formulate the complete convex programs with the objective and all the constraints needed to obtain the results of \secref{exem} as these problems are decoupled.

This allows to formulate the invariance of a convex set in terms of its support function
if it is differentiable.
We generalize this result with a relaxed notion of differentiability in \theoref{pdiff}.

\begin{mytheo}
  \label{theo:invexposed}
  Consider a nonempty closed convex set $\Csetvar$ such that
  $\supfun{\cdot}{\Csetvar}$ is differentiable.
  Then $\Csetvar$ is invariant for system \eqref{eq:as}
  with matrices $C, E \in \mathbb{R}^{r \times n}$
  if and only if
  \begin{equation}
    \label{eq:invgrad}
    \forall z \in \mathbb{R}^r, \la z, C \grad \supfun{E^\Tr z}{\Csetvar} \ra \le 0.
  \end{equation}
  \begin{proof}
    By \propref{exposed}, $F_\Csetvar(E^\Tr z) = \{\grad \supfun{E^\Tr z}{\Csetvar}\}$
    hence \eqref{eq:invface} is equivalent to \eqref{eq:invgrad}.
  \end{proof}
\end{mytheo}

\subsection{Ellipsoidal controlled invariant set}
\label{sec:ell}



In this section, we particularize \theoref{invexposed} to the case of ellipsoids.
Since the support function of an ellipsoid $\mathcal{E}_P$ is $\supfun{y}{\mathcal{E}_P} = \sqrt{y^\Tr P^{-1} y}$,
we have the following corollary of \theoref{invexposed}.
\begin{mycoro}[\cite{barmish1985necessary}]
  Given a positive definite matrix $P$, the ellipsoid $\mathcal{E}_P$ is
  controlled invariant for system~\eqref{eq:as} if and only if
  \begin{equation}
    \label{eq:invLMI}
    CP^{-1}E^\Tr + EP^{-1}C^\Tr \preceq 0.
  \end{equation}
\end{mycoro}

Observe that for the trivial case $\Image(B) = \mathbb{R}^n$ for system~\eqref{eq:cs},
\propref{csas} produces a system~\eqref{eq:as} with $r = 0$ hence the
LMI~\eqref{eq:invLMI} will be trivially satisfied for any $P^{-1}$, which is expected.

In comparison to \eqref{eq:controlled_LMI}, for a system~\eqref{eq:as} with matrices
$C, E \in \mathbb{R}^{r \times n}$, the LMI~\eqref{eq:controlled_LMI} has
size $n \times n$ while the LMI~\eqref{eq:invLMI} has only size $r \times r$.
The characterization of controlled invariance of ellipsoids using \eqref{eq:invLMI}
can also be obtained by applying an elimination procedure to reduce \eqref{eq:controlled_LMI};
see \cite[Equation~(7.11)]{boyd1994linear}.
However, uncertain or switched system may need a nonlinear state feedback to be quadratically
stabilizabilized~\cite{petersen1985quadratic}.
For such systems, \eqref{eq:controlled_LMI} is conservative since it assumes a linear feedback
while \eqref{eq:invLMI} does not assume anything about the feedback.
It was shown in \cite{barmish1985necessary} that if \eqref{eq:invLMI} is satisfied then a stabilizing nonlinear continuous
state feedback can be deduced from the solution $P$.
There is even a closed form for the feedback in case of single input~\cite[Eq.~(15)]{barmish1985necessary}.

\subsection{Polynomial controlled invariant set}
\label{sec:poly}
In this section, we derive the algebraic condition for
the controlled invariance of a set with polynomial support function.
This template is referred to as \emph{polyset}; see \cite[Section~1.5.3]{legat2020set}.
\begin{mycoro}
  \label{coro:poly}
  Given a homogeneous\footnote{A polynomial is \emph{homogeneous} if all its monomials have the same total degree} nonnegative polynomial $p(x)$ of degree $2d$, the set $\Csetvar$
  defined by the support function
  $\supfun{y}{\Csetvar} = p(y)^{\frac{1}{2d}}$
  is invariant for system~\eqref{eq:as}
  with matrices $C, E \in \mathbb{R}^{r \times n}$
  if and only if the polynomial
  \begin{equation}
    \label{eq:invSOS}
    z^\Tr C \grad p(E^\Tr z)
  \end{equation}
  is nonpositive for all $z \in \mathbb{R}^r$.
  \begin{proof}
    We have
    \[
      \grad \supfun{y}{\Csetvar} = \frac{1}{p(y)^{1-\frac{1}{2d}}} \grad p(y).
    \]
    If $p(y)$ is identically zero, this is trivially satisfied.
    Otherwise, $p(y)^{1-\frac{1}{2d}}$ is nonnegative and is zero
    in an algebraic variety of dimension $n - 1$ at most.
    Therefore,
    \eqref{eq:invgrad} is equivalent to \eqref{eq:invSOS}.
  \end{proof}
\end{mycoro}

While verifying the nonnegativity of a polynomial is co-NP-hard,
a sufficient condition can be obtained via the standard Sum-of-Squares programming framework; see ~\cite{blekherman2012semidefinite}.

\subsection{Piecewise semi-ellipsoidal controlled invariant set}
\label{sec:pell}

\cite{johansson1998computation} study the computation of piecewise quadratic Lyapunov functions for continuous-time autonomous piecewise affine systems.
\cite{legat2020piecewise} present a convex programming approach
to compute \emph{piecewise semi-ellipsoidal} controlled invariant sets of discrete-time control systems.
In this section, we show that \theoref{invexposed} provides the corresponding condition for continuous-time.

We recall \cite[Definition~2]{legat2020piecewise} below.
\begin{mydef}
  A \emph{polyhedral conic partition} of $\R^n$ is a set of $m$ polyhedral cones $(\mathcal{P}_i)_{i=1}^m$
  with nonempty interior
  such that for all $i \neq j$, $\dim(\mathcal{P}_i \cap \mathcal{P}_j) < n$
  and $\cup_{i=1}^m \mathcal{P}_i = \R^n$.
\end{mydef}

Piecewise semi-ellipsoids have a support function of the form
\begin{equation}
  \supfun{y}{\Csetvar} = \sqrt{y^\Tr Q_i y} \qquad y \in \mathcal{P}_i
\end{equation}
where $(\mathcal{P}_i)_{i=1}^m$ is a polyhedral conic partition.
The support function additionally has to satisfy \cite[(2) and (3)]{legat2020piecewise}
to ensure its continuity and convexity.

\begin{mytheo}
  \label{theo:pdiff}
  Consider a polyhedral conic partition $(\mathcal{P}_i)_{i=1}^m$ and
  a nonempty closed convex set $\Csetvar$
  such that 
  \[ \supfun{y}{\Csetvar} = f_i(y) \qquad y \in \mathcal{P}_i \]
  for differentiable functions $f_i : \mathcal{P}_i \to \R$. 
  The set $\Csetvar$ is invariant for system~\eqref{eq:as}
  with matrices $C, E \in \mathbb{R}^{r \times n}$
  if and only if, for all $i = 1, \ldots, m$ and $z \in \mathbb{R}^r$
  such that $E^\Tr z \in \mathcal{P}_i$, 
  we have
  \begin{equation}
    \label{eq:invpgrad}
    \la z, C \grad f_i(E^\Tr z) \ra \le 0.
  \end{equation}
  \begin{proof}
    Given $z \in \mathbb{R}^r$ such that $E^\Tr z$ is in the intersection of the boundary of $\Csetvar$ and
    the interior of $\mathcal{P}_i$,
    the support function is differentiable at $E^\Tr z$
    hence, by \propref{exposed},
    $F_\Csetvar(E^\Tr z) = \{\grad f_i(E^\Tr z)\}$.
    The condition~\eqref{eq:invface} is therefore reformulated as \eqref{eq:invpgrad}.

    Given a subset $I$ of $\{1, \ldots, m\}$ and
    $z \in \mathbb{R}^r$ such that $E^\Tr z$ is in
    the intersection of the boundary of $\Csetvar$ and $\cap_{i \in I} \mathcal{P}_i$,
    $F_{\Csetvar}(E^\Tr z)$ is the convex hull of
    $\grad \supfun{E^\Tr z}{\Csetvar}$ for each $i \in I$.
    For any convex combination (i.e., nonnegative numbers summing to 1) $(\lambda_i)_{i \in I}$,
    \eqref{eq:invpgrad} implies that
    \[ \la z, C \sum_{i \in I} \lambda_i \grad f_i(E^\Tr z) \ra = \sum_{i \in I} \lambda_i \la z, C \grad f_i(E^\Tr z) \ra \le 0. \]
  \end{proof}
\end{mytheo}

\begin{mycoro}
  A piecewise semi-ellipsoid $\Csetvar$
  is invariant for system~\eqref{eq:as}
  with matrices $C, E \in \mathbb{R}^{r \times n}$
  if and only if the quadratic form
  \begin{equation}
    \label{eq:invpell}
    z^\Tr CP_i^{-1}E^\Tr z + z^\Tr EP_i^{-1}C^\Tr z
  \end{equation}
  is nonpositive for all
  $i = 1, \ldots, m$ and
  $z \in \mathbb{R}^r$ such that $E^\Tr z \in \mathcal{P}_i$.
\end{mycoro}

The condition \eqref{eq:invpell} amounts to verifying the positive semidefiniteness of a quadratic form
when restricted to a polyhedral cone. When this cone is the positive orthant, this is called the
\emph{copositivity} which is co-NP-complete to decide \pcite{murty1987some}.
However, a sufficient LMI is given in \cite[Proposition~2]{legat2020piecewise}
and a necessary and sufficient condition is given by a hierarchy of Sum-of-Squares programs \cite[Chapter~5]{parrilo2000structured}.

\section{Numerical example}
\label{sec:exem}

\begin{figure}[!ht]
  \centering
  \polyset{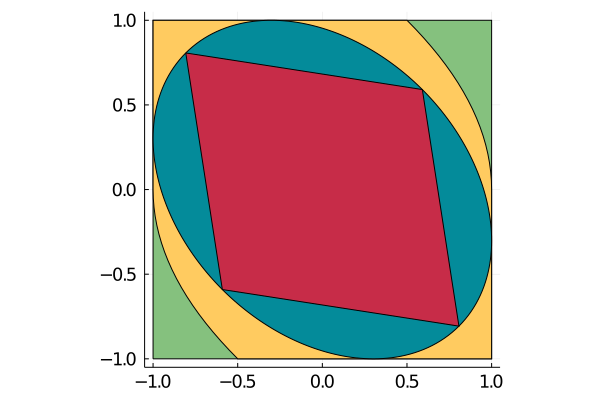}
  \polyset{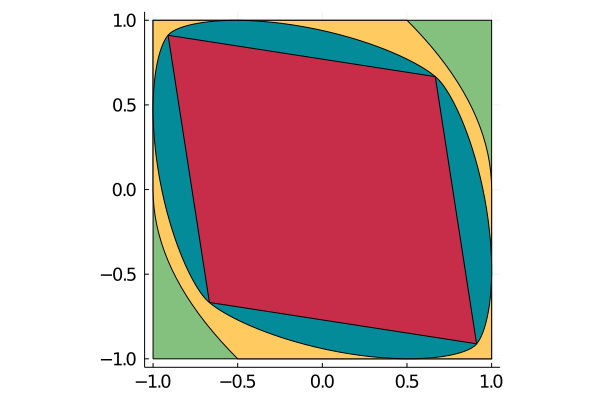}
  \polyset{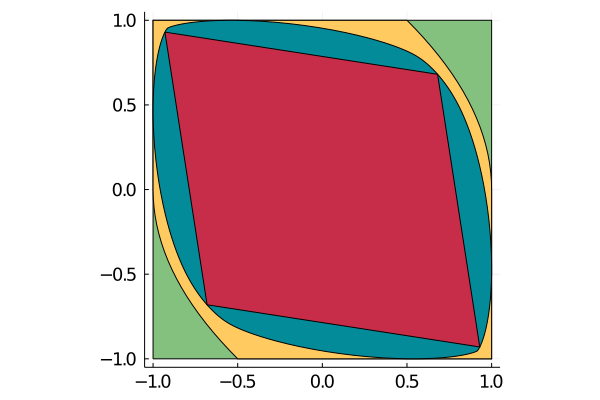}
  \polyset{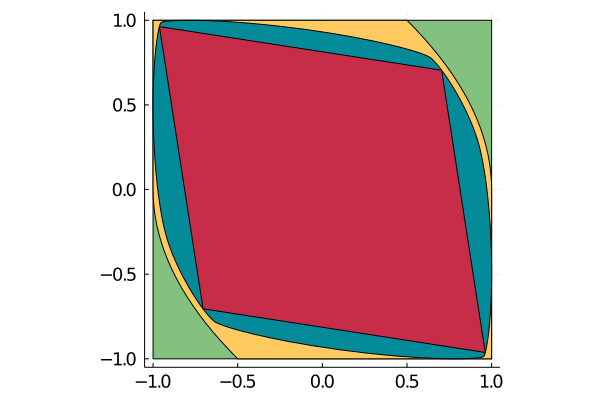}
  \polyset{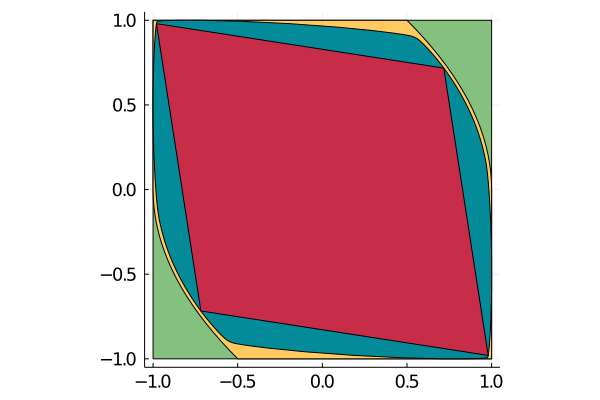}
  \caption{
    In blue are the solution for polysets of different degrees.
    The degrees from top to bottom are respectively 2, 4, 6, 10 and 20.
    The green set is the safe set $[-1, 1]^2$, the yellow set is the maximal
    controlled invariant set and the red set is $\gamma \mathcal{D}$.
    The sets are represented in the primal space in left figures
    and in polar space in the right figures.
  }
  \label{fig:poly}
\end{figure}

\begin{figure}[!ht]
  \centering
  \polyset{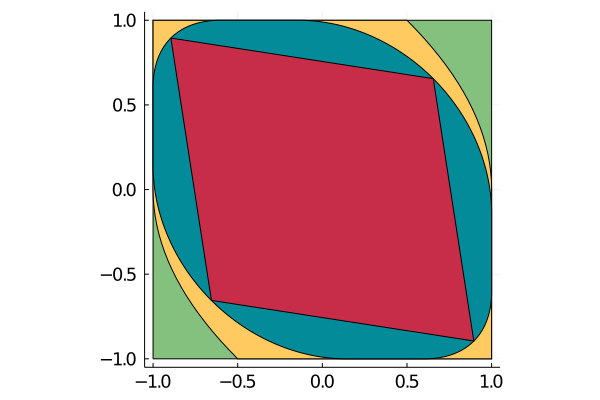}
  \polyset{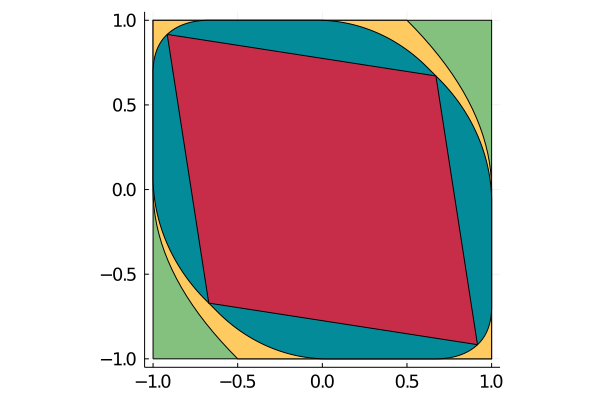}
  \caption{
    In blue are the solution for piecewise semi-ellipsoids for two
    different polyhedral conic partitions.
    The partitions from top to bottom are as described in \eqref{eq:sphere}
    with $m = (4, 3)$ (resp. $(8, 5)$).
    The green set is the safe set $[-1, 1]^2$, the yellow set is the maximal
    controlled invariant set and the red set is $\gamma \mathcal{D}$.
    The sets are represented in the primal space in left figures
    and in polar space in the right figures.
  }
  \label{fig:pell}
\end{figure}

In this section, we study a simple numerical example to illustrate our new approach.
Suppose we want to compute a controlled invariant set for the control system
\[
  \dot{x} =
  \begin{pmatrix}
    0 & 1\\
    0 & 0\\
  \end{pmatrix}
  x +
  \begin{pmatrix}
    0\\
    1
  \end{pmatrix}
  u
\]
with the state constraint $x \in [-1, 1]^2$ and input constraint $u \in [-1, 1]$.
We represent the state set $[-1, 1]^2$ and its polar in green in \figref{poly} and \figref{pell}.

The union of controlled invariant sets is controlled invariant.
Moreover, by linearity, the convex hull of the unions of controlled
invariant sets is controlled invariant.
Therefore, there exists a \emph{maximal} controlled invariant, i.e., a controlled invariant set in which all controlled invariant sets are included, for any family that is closed under union (resp. convex hull);
it is the union (resp. convex hull) of all controlled invariant sets included in $[-1, 1]^2$.

For this simple planar system, the maximal controlled invariant set
can be obtained by hand, it is 
\begin{equation*}
  \setbuild{(x_1, x_2) \in [-1, 1]^2}{x_1x_2 \le 0 \text{ or } |x_1| \le 1 - x_2^2/2}
\end{equation*}
Its polar is given by
\begin{align*}
  & \,\setbuild{x}{x_1x_2 \le 0 \text{ and } |x_1 - x_2| \le 1}\\
  \cup & \,\setbuild{x}{x_1(x_1 - x_2) \le 0 \text{ and } |x_1/2 + x_2| \le 1}\\
  \cup & \,\setbuild{x}{x_2(x_2 - x_1) \le 0 \text{ and } (2x_1 - \sign(x_1))^2 + 2x_2^2 \le 1}
\end{align*}
We represent it in yellow in \figref{poly} and \figref{pell}.

As \propref{csas} requires the input to be unconstrained,
we will consider projections onto the first two dimensions
of controlled invariant sets of the following control system:
\[
  \dot{x} =
  \begin{pmatrix}
    0 & 1 & 0\\
    0 & 0 & 1\\
    0 & 0 & 0
  \end{pmatrix}
  x +
  \begin{pmatrix}
    0\\
    0\\
    1
  \end{pmatrix}
  u.
\]
with the state constraint $x \in [-1, 1]^3$
and no input constraint.

Following \propref{csas}, we consider the algebraic system
\[
  \begin{pmatrix}
    1 & 0 & 0\\
    0 & 1 & 0\\
  \end{pmatrix}
  \dot{x} =
  \begin{pmatrix}
    0 & 1 & 0\\
    0 & 0 & 1\\
  \end{pmatrix}
  x
\]
with the state constraint $x \in [-1, 1]^3$.

While the \emph{maximal} invariant set is well defined,
it is not the case anymore when we restrict the set to belong to the family of ellipsoids,
polysets or piecewise semi-ellipsoids for a fixed polyhedral conic partition
as these families are not invariant under union nor convex hull.
The objective used to determine which invariant set is selected
depends on the particular application.
Let $\mathcal{D}$ be the convex hull of $\{(-1 + \sqrt3, -1 + \sqrt3), (-1, 1), (1 - \sqrt3, 1 - \sqrt3), (1, -1)\}$.
For this example, we maximize $\gamma$ such that $\gamma\mathcal{D}$ is
included in the projection of the invariant set onto the first two dimensions.
We represent $\gamma\mathcal{D}$ in red in \figref{poly} and \figref{pell}.

For the ellipsoidal template considered in \secref{ell}, the optimal solution
is shown in \figref{poly} as ellipsoids corresponds to polysets of degree 2.
The optimal objective value is $\gamma \approx 0.81$.

For the polyset template considered in \secref{poly}, the optimal solution
are represented in \figref{poly}. 
The optimal objective value for degree 4 (resp. 6, 10 and 20)
is $\gamma \approx 0.91$.
(resp.
$\gamma \approx 0.93$,
$\gamma \approx 0.96$ and
$\gamma \approx 0.98$).

For the piecewise semi-ellipsoidal template, we consider
polyhedral conic partitions
made of the conic hull of each facet of the polytope with extreme points
\begin{equation}
  \label{eq:sphere}
  (\cos(\alpha)\cos(\beta), \sin(\alpha)\cos(\beta), \sin(\beta))
\end{equation}
where $\alpha = 0, 2\pi/m_1, 4\pi/m_1, \ldots, 2(m_1-1)\pi/m_1$ and
$\beta = -\pi/2, \ldots, -2\pi/(m_2-1), -\pi/(m_2-1), 0, \pi/(m_2-1), 2\pi/(m_2-1), \ldots, \pi/2$.

The optimal objective value for $m = (4, 3)$ (resp. $(8, 5)$) is $\gamma \approx 0.89$ (resp. $\gamma \approx 0.92$).
The corresponding optimal solution is shown in \figref{pell}.

\section{Conclusion}
We proved a condition for continuous-time controlled invariance of a convex set
based on its support function.
We particularized the condition for three templates:
ellipsoids, polysets and piecewise semi-ellipsoids.
In the ellipsoidal case, it reduces to a known LMI,
in the polyset case, it provides a condition significantly less conservative than the existing one\footnote{Indeed, it is equivalent to invariance by \cororef{poly} and we showed in \secref{alg} that the existing approach is quite conservative.}
and in the piecewise semi-ellipsoidal case, it provides the first convex programming approach for continuous-time controlled invariance to the best of our knowledge.

As future work, we aim to apply this framework to other families such as
the \emph{piecewise polysets} defined in \cite{legat2020set}.
Moreover, instead of considering a uniform discretization of the hypersphere
as in \eqref{eq:sphere}, a more adaptive methods could be considered.
The sensitivity information provided by the dual solution of the optimization program could for instance determine
which pieces of the partition should be refined.

Finally, as the discrete-time version of this work developed in \cite{legat2020sum, legat2020piecewise}
also requires the set to be represented by its support function for the
optimization program to be convex,
both these methods and the method of this paper could be combined to compute controlled invariant sets for
hybrid automata using the condition of this paper for the invariance subject to the dynamics of each mode
and the condition of \cite{legat2020sum, legat2020piecewise} for the invariance subject to each reset map.



\bibliography{biblio}

\appendix
\section{Block matrices}

We have the following result for block matrices.
\begin{myprop}
  \label{prop:block}
  Consider a symmetric matrix $A \in \R^{n \times n}$ and a matrix $B \in \R^{n \times m}$.
  If the matrix
  \[
    C = \begin{bmatrix}
      A & B\\
      B^\Tr & 0
    \end{bmatrix}
  \]
  is positive semidefinite then $B$ is zero.
  \begin{proof}
    If $C$ is positive semidefinite,
    then there exists an integer $r$ and matrices $X \in \R^{r \times n}, Y \in \R^{r \times m}$
    such that
    \[
      C =
      \begin{bmatrix}
        X^\Tr\\
        Y^\Tr
      \end{bmatrix}
      \begin{bmatrix}
        X &
        Y
      \end{bmatrix}.
    \]
    From $Y^\Tr Y = 0$, we deduce that $Y = 0$ hence $B = X^\Tr Y = 0$.
  \end{proof}
\end{myprop}

\section{Convex geometry}
\label{sec:convex}

\begin{mydef}[{Support function \cite[p.~28]{rockafellar2015convex}}]
  Consider a convex set $\Csetvar$.
  The \emph{support function} of $\Csetvar$ is defined as
  \[ \supfun{y}{\Csetvar} = \sup_{x \in \Csetvar} \la y, x \ra. \]
\end{mydef}

\begin{mydef}[{Polar set}]
  For any convex set $\Csetvar$ the polar of $\Csetvar$,
  denoted $\polar{\Csetvar}$,
  is defined as
  \[
    \polar{\Csetvar} = \setbuild{y}{\supfun{y}{\Csetvar} \le 1}.
  \]
\end{mydef}

We define the \emph{tangent cone} as follows \cite[Definition~4.6]{blanchini2015set}.

\begin{mydef}[{Tangent cone}]
  Given a closed convex set $\setS$ and a distance function $\dist{\setS}{x}$,
  the \emph{tangent cone} to $\setS$ at $x$ is defined as follows:
  \[ T_{\setS}(x) = \left\{\, y \mid \lim_{\tau \to 0} \frac{\dist{\setS}{x + \tau y}}{\tau} = 0 \,\right\} \]
  where the distance is defined as
  \[
    \dist{\setS}{x} = \inf_{y \in \setS} \|x - y\|
  \]
  where $\|\cdot\|$ is a norm.
  The tangent cone is a convex cone and is independent of the norm used.

\end{mydef}

For a convex set $\Csetvar$, the \emph{normal cone} is the polar of the tangent cone $N_\Csetvar(x) = \polar{T_\Csetvar}(x)$.

The \emph{exposed face} (also called the \emph{support set}, e.g., in \cite[Section~1.7.1]{schneider2013convex})
is defined as follows \cite[Definition~3.1.3]{hiriart2012fundamentals}.
\begin{mydef}[{Exposed face}]
  Consider a nonempty closed convex set $\Csetvar$.
  Given a vector $y \neq 0$, the \emph{exposed face} of $\Csetvar$ associated
  to $y$ is
  \[
    F_\Csetvar(y) = \{\, x \in \Csetvar \mid \la x, y \ra = \supfun{y}{\Csetvar} \,\}.
  \]
\end{mydef}

The exposed faces and normal cones are related by the following property \cite[Proposition~C.3.1.4]{hiriart2012fundamentals}.

\begin{myprop}
  \label{prop:normalsupfun}
  Consider a nonempty closed convex set $\Csetvar$.
  For any $x \in \Csetvar$ and nonzero vector $y$,
  $x \in F_\Csetvar(y)$ if and only if $y \in N_\Csetvar(x)$.
\end{myprop}

Given a set $\Sset$ and a matrix $A$,
let $A^{-\Tr}$ denote the preimage
$\setbuild{x}{A^\Tr x \in \Sset}$.

\begin{myprop}[{\cite[Corollary~16.3.2]{rockafellar2015convex}}]
  \label{prop:podu}
  For any convex set $\Csetvar$ and linear map $A$,
  \begin{align*}
    \polar{(A\Csetvar)} & = A^{-\Tr} \polar{\Csetvar}.
  \end{align*}
  where $\polar{\Csetvar}$ denotes the polar of the set $\Csetvar$.
\end{myprop}

When the support function is differentiable at a given point, $F_\Csetvar$ is a singleton and may be directly obtained using the following result:
\begin{myprop}[{\cite[Corollary~25.1.2]{rockafellar2015convex}}]
  \label{prop:exposed}
  Given a nonempty closed convex set $\Csetvar$,
  if $\supfun{y}{\Csetvar}$ is differentiable at $y$ then $F_{\Csetvar}(y) = \{\grad \supfun{y}{\Csetvar}\}$.
\end{myprop}

In fact,
for nonempty compact convex sets,
the differentiability at $y$ is even
a necessary and sufficient conditions for
the unicity of $F_{\Csetvar}(y)$~\cite[Corollary~1.7.3]{schneider2013convex}.

\end{document}